\documentclass[11pt,twoside]{amsart}

\usepackage{amsmath,amssymb,amsfonts}

\title[Tensor product of stably equivalent algebras]{On a question of Rickard on tensor product of stably equivalent algebras}

\author{Serge Bouc}
\address{\newline
S.B.: Universit\'e de Picardie,
\newline LAMFA (UMR 7352 du CNRS),
\newline 33 rue St Leu,
\newline F-80039 Amiens Cedex 1,
\newline France}
\email{serge.bouc@u-picardie.fr}

\author{Alexander Zimmermann}
\address{\newline
A.Z.: Universit\'e de Picardie,
\newline D\'epartement de Math\'ematiques et LAMFA (UMR 7352 du CNRS),
\newline 33 rue St Leu,
\newline F-80039 Amiens Cedex 1,
\newline France}
\email{alexander.zimmermann@u-picardie.fr}

\newtheorem*{Theo2}{{Theorem}}
\newtheorem{Lemma1}{{Lemma}}
\newtheorem{Theo1}[Lemma1]{{Theorem}}
\newtheorem{Def1}[Lemma1]{{Definition}}
\newtheorem{Prop1}[Lemma1]{{Proposition}}
\newtheorem{Claim1}[Lemma1]{{Claim}}
\newtheorem{Rem1}[Lemma1]{{Remark}}
\newtheorem{Cor1}[Lemma1]{{Corollary}}
\newtheorem{Ex1}[Lemma1]{{Example}}
\newtheorem{Qu1}[Lemma1]{{Question}}

\newenvironment{Lemma}{\begin{Lemma1}}{\end{Lemma1}}
\newenvironment{Def}{\begin{Def1}\em}{\end{Def1}}
\newenvironment{Prop}{\begin{Prop1}}{\end{Prop1}}
\newenvironment{Rem}{\begin{Rem1}\rm}{\end{Rem1}}
\newenvironment{Theorem}{\begin{Theo1}}{\end{Theo1}}

\newenvironment{Cor}{\begin{Cor1}}{\end{Cor1}}

\newcommand{\lra}{\longrightarrow}

\newcommand{\ra}{\rightarrow}
\newcommand{\sdp}{\times\kern-.2em\vrule height1.1ex depth-.05ex}
\newcommand{\epi}{\lra \kern-.8em\ra}

\newcommand{\N}{{\mathbb N}}

\newcommand{\F}{{\mathbb F}}

\newcommand{\Z}{{\mathbb Z}}

 \normalbaselines

\date{\today}

\thanks{Both authors were supported by a grant STIC Asie 'Escap' from the Minist\`ere des Affaires \'Etrang\`eres de la France}

\begin{document}

\begin{abstract} Let $\overline\F_p$ be the algebraic closure of the prime field of characteristic $p$.
After observing that the principal block $B$ of $\overline\F_pPSU(3,p^r)$ 
is stably equivalent of Morita type to its Brauer correspondent $b$, we 
show however that the centre of $B$ is not isomorphic as an algebra to 
the centre of $b$ in the cases $p^r\in\{3,4,5,8\}$. As a consequence, 
the algebra $B\otimes_{\overline{\F}_p}\overline \F_p[X]/X^p$ is not 
stably equivalent of Morita type to $b\otimes_{\overline\F_p}\overline\F_p[X]/X^p$ in these cases. 
This yields a negative answer to a question of Rickard.
\end{abstract}

\maketitle

\section*{Introduction}

Let $K$ be a field, and let $A$, $B$, $C$ and $D$ be finite dimensional $K$-algebras.
Rickard showed in \cite{derivedfunctors} that if $A$ and $B$ are derived
equivalent, and if $C$ and $D$ are derived equivalent, then also $A\otimes_KC$
and $B\otimes_KD$ are derived equivalent. Rickard asks in \cite[Question 3.8]{recentadvances} if this still holds when replacing
derived equivalence by stable equivalence of Morita type. It is clear that we
have to suppose that all algebras involved have no semisimple direct factor. A
result due to Liu \cite{summands} shows that then we may suppose that all
algebras are indecomposable. In \cite{stabletensor} Liu, Zhou and the second
author showed that the question has a negative solution in case $A$, $B$, $C$
and $D$ are not necessarily selfinjective. However, a derived equivalence
between selfinjective algebras $A$ and $B$ induces a stable equivalence
of Morita type between $A$ and $B$. If $A$ and $B$ are not selfinjective,
then this implication is not valid. Hence, the natural playground for Rickard's
question are selfinjective algebras.

The purpose of this paper is to give a counterexample to Rickard's question. For an algebraically closed base field $K$ of characteristic $p$ we construct
symmetric $K$-algebras $A$ and $B$ which are stably equivalent of Morita type,
but $A\otimes_KK[X]/X^p$ and $B\otimes_KK[X]/X^p$ are not stably equivalent of
Morita type.

Note that this answers the general case. Indeed,
if $A\otimes_K C$ is stably equivalent of Morita type to $B\otimes_KC$ and
$B\otimes_KC$ is stably equivalent of Morita type to $B\otimes_KD$ then $A\otimes_KC$ is stably equivalent of Morita type to $B\otimes_KD$. Hence,
we may suppose $C=D$.

In recent years many attempts were proposed to lift a stable
equivalence of Morita type between selfinjective algebras to a
derived equivalence. It is known that this is not possible in
general, as is seen by the mod $2$ group ring of a dihedral group
of order $8$ and the stable equivalence induced by a uniserial
endotrivial module of Loewy length $3$. This was used in
\cite{stabletensor} for example. In this paper we give a new
incidence of this fact. Moreover, we provide two symmetric
algebras, which are stably equivalent of Morita type, and
have non isomorphic centres.

Our example is the principal $p$-block of the group $PSU(3,p^r)$
and its Brauer
correspondent for $p^r\in\{3,4,5,8\}$.

We recall in the first section some basic facts and results which
we need for
our construction. In Section~\ref{Examplesection}
we give our main result and its proof,
and in Section~\ref{programssection} we display the GAP program
needed for the proof. In Section~\ref{Nsection} we determine the
algebraic structure
of the centre of $KN_G(S)$ for $G=PSU(3,p^r)$ and $S$ one of
its Sylow $p$-subgroups for all primes $p$ and integers $r$.

{\bf Acknowledgement:} The idea of this paper was born during a visit
of both authors in Beijing Normal University. We are very grateful to
Yuming Liu for his great hospitality. Moreover, we thank Yuming Liu
for suggesting a question to us leading to the present paper, and also for pointing out that it is not sufficient to show an abstract non isomorphism of the centres of the blocks. This led us to complete our proof by adding Lemma~\ref{nontrivialtensorproduct}.

\section{Background}

Recall the following

\begin{Def} \cite{Broue}, (cf also \cite[Chapter 5]{reptheo})
Let $A$ and $B$ be two finite dimensional algebras over a field $K$.
Then $A$ and $B$ are stably equivalent of Morita type if there is an  $A\otimes_KB^{op}$-module $M$ and a $B\otimes_KA^{op}$-module $N$ such that
\begin{itemize}
\item $M$ is projective as $A$-module, and as $B^{op}$-module
\item $N$ is projective as $A^{op}$-module and as $B$-module
\item there is a projective $A\otimes_K A^{op}$-module $P$ and a projective $B\otimes_K B^{op}$-module $Q$ such that $M\otimes_BN\simeq B\oplus Q$ as $B\otimes_KB^{op}$-modules and $N\otimes_AM\simeq A\oplus P$ as $A\otimes_KA^{op}$-modules.
\end{itemize}
\end{Def}

Independently Rickard \cite{Ristable} as well as Keller and Vossieck \cite{KV}, show that if $A$ and $B$ are derived equivalent selfinjective algebras, then $A$ and $B$ are stably equivalent of Morita type.

Brou\'e defined $Z^{st}(A):=\underline{\text{End}}_{A\otimes_KA^{op}}(A)$ and
$$Z^{pr}(A):=\ker(\text{End}_{A\otimes_KA^{op}}(A)\rightarrow \underline{\text{End}}_{A\otimes_KA^{op}}(A))$$
where we denote by $\underline{\text{End}}$ the endomorphisms taken in the stable module category.

The centre of an algebra is an invariant of a derived equivalence, as was shown by Rickard. The stable centre $Z^{st}(A)$ is an important invariant under stable equivalences of Morita type, as was shown by Brou\'e.

\begin{Prop} (Brou\'e~\cite{Broue}; see also \cite[Chapter 5]{reptheo})
If $A$ and $B$ are stably equivalent of Morita type, then $Z^{st}(A)\simeq Z^{st}(B)$ as algebras.
\end{Prop}

Now, Liu, Zhou and the second author give a criterion to determine the dimension of $Z^{st}(A)$.

\begin{Theorem} \cite[Proposition 2.3 and Corollary 2.7]{LZZ}
\label{rankCartanprojectivecentre}
Let $A$ be a finite dimensional symmetric algebra over an algebraically closed field $K$ of characteristic $p>0$. Then $\dim_K(Z^{pr}(A))=\text{\rm rank}_p(C_A)$ where $C_A$ is the Cartan matrix of $A$ and where $\text{\rm rank}_p(C_A)$ denotes its rank as matrix over $K$.
\end{Theorem}

Moreover, we recall a conjecture of Auslander-Reiten. In \cite{ARS} Auslander and Reiten conjecture that if $A$ and $B$ are stably equivalent finite dimensional algebras, then the number of simple non-projective $A$-modules and the number of non-projective simple $B$-modules coincides. Again in \cite{LZZ} we show

\begin{Theorem} \cite[Theorem 1.1]{LZZ}
Let $K$ be an algebraically closed field and let $A$ and $B$ be two finite dimensional $K$-algebras, which are stably equivalent of Morita type and which do not have any semisimple direct factor. Then the number of isomorphism classes of non-projective simple $A$-modules is equal to the number of non-projective simple $B$-modules if and only if $\dim_K(HH_0(A))=\dim_K(HH_0(B))$, where $HH_0$ denotes the degree $0$ Hochschild homology.
\end{Theorem}

In particular, if $A$ and $B$ are symmetric, then Hochschild homology and cohomology coincides, and the number of non-projective simple $A$-modules is equal to the number of non-projective simple $B$-modules if and only if the centres of $A$ and of $B$ have the same dimension.

The following lemma is well-known to the experts, but for the convenience of the reader, and since it is crucial to our arguments, we include the short proof. For an algebra $A$ denote by $J(A)$ its Jacobson radical.

\begin{Lemma}\label{radicaloftensorproduct}
Let $K$ be a perfect field and let $A$ and $B$ be finite dimensional $K$-algebras. Then $J(A\otimes_KB)=J(A)\otimes_KB+A\otimes_KJ(B)$.
\end{Lemma}

\begin{proof} 
It is clear that $J(A)\otimes_KB+A\otimes_KJ(B)$ is a nilpotent ideal of $A\otimes_KB$, and therefore we get $$J(A)\otimes_KB+A\otimes_KJ(B)\subseteq J(A\otimes_KB).$$

Now, $(A\otimes_KB)/(J(A)\otimes_KB+A\otimes_KJ(B))=A/J(A)\otimes_KB/J(B)$
and both $K$-algebras $A/J(A)$ and $B/J(B)$ are semisimple.  
Since $K$ is perfect, every finite extension $L$ of $K$ is a separable field extension.
By \cite[Corollary 7.6]{CR} a finite dimensional semisimple $K$-algebra $C$ is separable if and only if the centres of each of the Wedderburn components is a separable field extension of $K$. Hence $A/J(A)$ and $B/J(B)$ are both separable $K$-algebras. By \cite[Corollary 7.8]{CR} the algebra
$A/J(A)\otimes_KB/J(B)$ is semisimple. Therefore  
$$J(A)\otimes_KB+A\otimes_KJ(B)\supseteq J(A\otimes_KB).$$ 
This shows the statement.
\end{proof}

\begin{Rem} (cf e.g.~\cite[Example 1.7.17]{reptheo}) Lemma~\ref{radicaloftensorproduct} is wrong if we drop the assumption that $K$ is perfect: e.g. let $p$ be a prime, and $K=\F_p(U)$ be the field of rational fractions over the finite field $\F_p$. Let $A=K[X]/(X^p-U)$. Then $A$ is a purely inseparable extension of $K$, of dimension $p$. In particular it is a reduced (commutative) algebra, i.e. $J(A)=0$. But $A\otimes_KA\cong K[X,Y]/(X^p-U,Y^p-U)$ contains the non zero element $X-Y$, such that $(X-Y)^p=U-U=0$. Hence $J(A\otimes_KA)\neq 0$.
\end{Rem}

\begin{Lemma}\label{cartanmatrixoftensorproducts}
Let $K$ be an algebraically closed field, let $\Lambda$ and $\Delta$ be finite dimensional $K$-algebras, and suppose that $\Delta$ is local. Then the projective indecomposable $\Lambda\otimes_K\Delta$-modules are precisely the modules $P\otimes_K\Delta$ for
projective indecomposable $\Lambda$-modules $P$, and if $C_\Lambda$ is the Cartan matrix of
$\Lambda$, then the Cartan matrix of $\Lambda\otimes\Delta$ is $C_{\Lambda\otimes_K\Delta}=\dim_K(\Delta)\cdot C_\Lambda$.
\end{Lemma}

\begin{proof}
Let $P$ and $Q$ be a indecomposable projective $\Lambda$-modules.
Then $P\otimes_K\Delta$ is a projective indecomposable $\Lambda\otimes_K\Delta$-module. Indeed,
$\text{End}_{\Lambda\otimes_K\Delta}(P\otimes_K\Delta)\simeq \text{End}_{\Lambda}(P)\otimes_K\Delta^{op}$.
Moreover, since $\Gamma:={\rm End}_\Lambda(P)^{op}$ and $\Delta$ are local $K$-algebras their radical
quotient are finite-dimensional skew-fields, and therefore
$\Gamma/J(\Gamma)\simeq K\simeq\Delta/J(\Delta)$ since $K$ is algebraically closed. Moreover, by Lemma~\ref{radicaloftensorproduct} we get
$J(\Gamma\otimes_K\Delta)= J(\Gamma)\otimes \Delta+\Gamma\otimes_KJ(\Delta).$
On the other hand, $$(\Gamma\otimes_K\Delta)/\left(J(\Gamma)\otimes_K \Delta+\Gamma\otimes_KJ(\Delta)\right)=K\otimes_K K=K$$
and hence we get
$\Gamma\otimes_K\Delta$ is local, and therefore $P\otimes_K\Delta$ is indecomposable. Now,
$$\text{Hom}_{\Lambda\otimes_K\Delta}(P\otimes_K\Delta,Q\otimes_K\Delta)=
\text{Hom}_\Lambda(P,Q)\otimes_K\Delta^{op}.$$
Taking $K$-dimensions proves the lemma.
\end{proof}

\begin{Rem}\label{Zstofpolynomial}
As a special case of Lemma~\ref{cartanmatrixoftensorproducts} we get
$C_{A\otimes_KK[X]/X^p}=p\cdot C_A$ for algebraically closed fields
$K$ of characteristic $p$.
Hence we get by Theorem~\ref{rankCartanprojectivecentre} that
$Z^{pr}(A\otimes_KK[X]/X^p)=0$ for algebraically closed fields
$K$ of characteristic $p$ and symmetric $K$-algebras $A$.
\end{Rem}

\begin{Lemma}\label{nontrivialtensorproduct}
Let $K$ be a perfect field and let
$n,m$ be positive integers. Let $A$ and $B$ be finite
dimensional commutative $K$-algebras. If $J^{n+1}(A)=0\neq J^n(A)$
and $J^{m+1}(B)=0\neq J^m(B)$, then $$J^{n+m+1}(A\otimes_K B)=0\neq J^{n+m}(A\otimes_K B)=J^n(A)\otimes_K J^m(B).$$
\end{Lemma}

\begin{proof} By Lemma~\ref{radicaloftensorproduct}, we have $J(A\otimes_K B)=J(A)\otimes_K B+A\otimes_K J(B)$.
Therefore
$$J^{n+m+1}(A\otimes_K B)=\sum_{k=0}^{n+m+1}J^k(A)\otimes_K J^{n+m+1-k}(B)=0\;\;\;.$$
Similarly
$$J^{n+m}(A\otimes_K B)=\sum_{k=0}^{n+m}J^k(A)\otimes_K J^{n+m-k}(B)=J^n(A)\otimes_K J^m(B)\neq 0\;\;\;,$$
which completes the proof.
\end{proof}

\begin{Rem} \label{centreofAtensorfruncatedpoly}
Let $K$ be any field, and $A$ be a $K$-algebra. We give an elementary argument to determine the centre of $A\otimes_KK[X]/X^p$.
It is clear that $A\otimes_KK[X]/X^p\cong A[X]/X^p$. Now, let
$a:=a_0+a_1X+\dots a_{p-1}X^{p-1}\in A[X]$.
Then for $b:=b_0\in A\cdot 1$ we get
$$ab-ba=(a_0b-ba_0)+\dots+(a_{p-1}b-ba_{p-1})X^{p-1}$$
and so $a\in Z(A)$ implies that $a$ commutes with any $b\in A$, and hence $a_0,\dots,a_{p-1}$ are all in $Z(A)$.
Conversely, it is clear that $Z(A)[X]/X^p\subseteq Z(A[X]/X^p)$
since $aX^n$ commutes with all elements of $A[X]/X^p$ whenever $a\in A$ and
since sums of elements in the centre are still central.
\end{Rem}

\begin{Lemma}\label{lastloewyofpolynomials}
If $K$ is a perfect field and $A$ is a finite dimensional 
$K$-algebra, and if moreover
$J^n(Z(A))\neq 0=J^{n+1}(Z(A))$, then
$$0\neq J^{n+p-1}\big(Z(A\otimes_KK[X]/X^p)\big)=J^n\big(Z(A)\big)
\otimes_KX^{p-1}K[X]/X^{p}$$
and $$J^{n+p}\big(Z(A\otimes_KK[X]/X^p)\big)=0.$$
\end{Lemma}

\begin{proof} This is an immediate consequence of
Lemma~\ref{nontrivialtensorproduct}.
\end{proof}

\begin{Cor}\label{whatweshow}
Let $K$ be an algebraically closed field of characteristic $p>0$ and let $A$ and $B$ be two finite dimensional $K$-algebras and let $n,m\in\N$ such that
$J^n(Z(A))\neq 0=J^{n+1}(Z(A))$ and $J^m(Z(B))\neq 0=J^{m+1}(Z(B))$.
If $\dim_K(J^n(Z(A)))\neq\dim_K(J^m(Z(B))$, then $A\otimes_KK[X]/X^p$ and
$B\otimes_KK[X]/X^p$ are not
stably equivalent of Morita type.
\end{Cor}

\begin{proof}
If $n\neq m$, then $Z(A\otimes_KK[X]/X^p)\not\simeq Z(B\otimes_KK[X]/X^p)$
by Lemma~\ref{lastloewyofpolynomials} since the Loewy lengths of the centres
are different. If $n=m$, then
Lemma~\ref{lastloewyofpolynomials} shows that
the centres of $A\otimes_KK[X]/X^p$ and of $B\otimes_KK[X]/X^p$ are not
isomorphic since the dimension of the lowest Loewy layers of
the centres are not of the same dimension. 
Remark~\ref{Zstofpolynomial} shows that
$Z(A\otimes_KK[X]/X^p)=Z^{st}(A\otimes_KK[X]/X^p)$ and
$Z(B\otimes_KK[X]/X^p)=Z^{st}(B\otimes_KK[X]/X^p)$.
Since the stable centre is invariant under stable equivalence of Morita type, we get the statement.
\end{proof}

\begin{Rem} For a field $K$ and a $K$-algebra $A$ let $n_A$ be the number of
isomorphism classes of simple nonprojective $A$-modules.
Auslander-Reiten conjecture \cite[page 409, Conjecture (5)]{ARS} that if $A$
and $B$ are stably equivalent finite dimensional $K$-algebras, then $n_A=n_B$.
\cite[Theorem 1.1]{LZZ} shows that if $K$ is algebraically closed
and if $A$ and $B$ are indecomposable finite dimensional $K$-algebras which are stably equivalent of Morita type, then $n_A=n_B$ is equivalent to
$\dim_K(HH_0(A))=\dim_K(HH_0(B))$. If $A$ is symmetric, then there is a vector
space isomorphism $HH_0(A)\simeq HH^0(A)=Z(A)$, we
see that the
Auslander-Reiten conjecture implies that
$\dim_K\big(Z(A)\big)=\dim_K\big(Z(B)\big)$.
More precisely by \cite[Corollary 1.2]{LZZ}, for two indecomposable
symmetric algebras $A$ and $B$ over an algebraically closed field
$K$ we have $n_A=n_B\Leftrightarrow \dim_K\big(Z^{pr}(A)\big)=\dim_K\big(Z^{pr}(B)\big)$,
where by definition $Z^{st}(A)=Z(A)/Z^{pr}(A)$. The link to our proof is
now given by the fact that for every algebra the Higman ideal $H(A)$ of $A$ equal $Z^{pr}(A)$, and for symmetric algebras $A$
over an algebraically closed field $K$ we have
$\dim_K(H(A))$ equals the $p$-rank of $C_A$.
\end{Rem}

\section{The Example}
\label{Examplesection}

Let $\overline\F_p$ be the algebraic closure of the prime field $\F_p$ of characteristic $p$.
Let $q=p^n$ for some integer $n$.

We recall some results on the geometry of $PSU(n,q)$ (cf e.g. \cite[II Satz 10.12, page 242]{Huppert}).
The group $G:=PSU(3,q)$ acts doubly transitively on the unitary quadric $Q$ of cardinal $q^3+1$. Note that we use the GAP notation, not the notation used in \cite[II Satz 10.12, page 242]{Huppert}, namely, $PSU(3,q)$ is defined over a field with $q^2$ elements, and is a natural quotient of a subgroup of $SL_2(q^2)$ (and not of $SL_2(q)$ !).
The stabiliser of a point $X$ of $Q$ is the normaliser in $G$ of a Sylow $p$-subgroup $P$
of $G$. Therefore two different conjugate Sylow $p$-subgroups $P$ and $^gP$ of $G$
fix two different points $X$ and $gX$ of $Q$. Hence $^gP\cap P=1$ if
$g\not\in N_G(P)$, or in other words, $G$ has a trivial intersection Sylow $p$-subgroup structure.
%For the statement one may also refer to Bender~\cite{Bender}.
This implies that Green correspondence  gives a
stable equivalence of Morita type between the principal block $B$ of
$\overline{\F}_pG$ and its Brauer correspondent $b$
(cf e.g. \cite[Chapter 2, Proposition 2.1.23 and Proposition 2.4.3]{reptheo}).

The GAP
\cite{Gapprograms} program in Section 3 computes the Loewy series of
the ring $Z(\F_2PSU(3,4))$ and of $Z(\F_2N_{PSU(3,4)}(S))$ for some
Sylow $2$-subgroup of $PSU(3,4)$. Observe moreover that $\overline\F_2PSU(3,4)$ has two blocks, the principal one and another block of defect $0$ (corresponding to the Steinberg character). Moreover,
the dimensions of the Loewy series obtained over $\F_2$ also hold by
extending the scalars to $\overline\F_2$, using Lemma~\ref{radicaloftensorproduct}.

We obtain that
$$\dim_{\overline\F_2}(Z(B))=21=\dim_{\overline\F_2}(Z(b))$$
$$\dim_{\overline\F_2}(J(Z(B)))=20=\dim_{\overline\F_2}(J(Z(b)))$$
$$\dim_{\overline\F_2}(J^2(Z(B)))=5\neq 4=\dim_{\overline\F_2}(J^2(Z(b)))$$
$$\dim_{\overline\F_2}(J^3(Z(B)))=0=\dim_{\overline\F_2}(J^3(Z(b))).$$

Similarly we get for the centre of the principal block $B$ of $PSU(3,8)$ and the  centre of its Brauer correspondent $b$
$$\dim_{\overline\F_2}(Z(B))=27=\dim_{\overline\F_2}(Z(b))$$
$$\dim_{\overline\F_2}(J(Z(B)))=26=\dim_{\overline\F_2}(J(Z(b)))$$
$$\dim_{\overline\F_2}(J^2(Z(B)))=3\neq 2=\dim_{\overline\F_2}(J^2(Z(b)))$$
$$\dim_{\overline\F_2}(J^3(Z(B)))=0=\dim_{\overline\F_2}(J^3(Z(b))).$$

An immediate variant of the program shows that this is a quite
general phenomenon in odd characteristic.
The group $PSU(3,3)$ gives an example in characteristic $3$ since,
denoting by $B$ the principal block of ${\overline\F}_3PSU(3,3)$ and by
$b$ its Brauer correspondent,
$$\dim_{\overline\F_3}(Z(B))=13=\dim_{\overline\F_3}(Z(b))$$
$$\dim_{\overline\F_3}(J(Z(B)))=12=\dim_{\overline\F_3}(J(Z(b)))$$
$$\dim_{\overline\F_3}(J^2(Z(B)))=4\neq 3=\dim_{\overline\F_3}(J^2(Z(b)))$$
$$\dim_{\overline\F_3}(J^3(Z(B)))=0=\dim_{\overline\F_3}(J^3(Z(b))).$$
The group $PSU(3,5)$ gives an example in characteristic $5$ since,
denoting by $B$ the principal block of ${\F}_5PSU(3,5)$ and by
$b$ its Brauer correspondent,
$$\dim_{\overline\F_5}(Z(B))=13=\dim_{\overline\F_5}(Z(b))$$
$$\dim_{\overline\F_5}(J(Z(B)))=12=\dim_{\overline\F_5}(J(Z(b)))$$
$$\dim_{\overline\F_5}(J^2(Z(B)))=2\neq 1=\dim_{\overline\F_5}(J^2(Z(b)))$$
$$\dim_{\overline\F_5}(J^3(Z(B)))=0=\dim_{\overline\F_5}(J^3(Z(b))).$$

\begin{Theorem}\label{notstablyequivalent}
Let $K$ be the algebraic closure of $\F_p$ and let $B$ be the principal
block of $PSU(3,p^r)$. Let $b$ be the Brauer correspondent of $B$ in the
group ring of the normaliser of a $2$-Sylow subgroup of $PSU(3,p^r)$. Then
$B$ and $b$ are stably equivalent of Morita type. If moreover
$p^r\in\{3,4,5,8\}$, then the square of the
Jacobson radical of $Z(B)$ is of different dimension than the
square of the Jacobson radical of  $Z(b)$, whereas $Z(B)$ and
$Z(b)$ both have Loewy length $3$. In
particular $B\otimes_KK[X]/X^p$ is not stably equivalent of Morita
type to $b\otimes_KK[X]/X^p$.
\end{Theorem}

\begin{proof}
As seen at the beginning of this section, $B$ and $b$ are stably equivalent of Morita type
by Green correspondence.

The GAP \cite{Gapprograms} program in Section~\ref{programssection}
shows that the Loewy series of the centres of $B$ and of $b$ are of the
same length but the dimensions of the Loewy layers are not equal.
In particular the lowest Loewy layers of the algebras $Z(B)$ and $Z(b)$
have different dimension.

Corollary~\ref{whatweshow}
implies that $B\otimes_{K}K[X]/X^p$ is not stably equivalent of Morita type to $b\otimes_{K}K[X]/X^p$.
\end{proof}

\newcommand{\fracb}[2]{\frac{\displaystyle #1}{\displaystyle #2}}
\begin{Rem} The above examples suggest that in general, with the notation of Theorem~\ref{notstablyequivalent}, the dimension of $J^2\big(Z(B)\big)$ could always be equal to $1+\dim_KJ^2(Z(b))$. By Theorem~\ref{J2ZkN}, this is equal to $\fracb{p^r+1}{\gamma\rule{0ex}{1.5ex}}$, where $\gamma$ is the greatest common divisor of $p^r+1$ and $3$.
\end{Rem}

\section{The GAP program}
\label{programssection}

We display here the GAP program we used.

\bigskip
\begin{footnotesize}
\begin{verbatim}
# the characteristic p
prem:=2;
#
# The group G
g:=PSU(3,prem^2);
#
# the ground field k
corps:=GF(prem);
#
s:=SylowSubgroup(g,prem);
# the normalizer NS of a Sylow p-subgroup
ns:=Normalizer(g,s);
#
# getting a permutation representation of G of smaller degree
f:=FactorCosetAction(g,ns);
g:=Image(f);
ns:=Image(f,ns);
#
# uncomment next line to replace G by NS
#g:=ns;
#
# computing the structure constants of ZkG
c:=ConjugacyClasses(g);
rc:=List(c,Representative);
lc:=Length(c);
ci:=List([1..lc],x->First([1..lc],y->rc[x]^(-1) in c[y]));
l:=List([1..lc],x->NullMat(lc,lc,corps));
for iu in [1..lc] do
    u:=c[iu];
    if rc[iu]=One(g) then
        for iv in [iu..lc] do
            Print("\r",iu,":",iv,"/",lc,"   ");
            v:=List([1..lc],x->Zero(corps));
            v[iv]:=One(corps);
            l[iu][iv]:=v;
            l[iv][iu]:=v;
        od;
    else
        for iv in [iu..lc] do
            Print("\r",iu,":",iv,"/",lc,"   ");
            w:=c[ci[iv]];
            v:=List(List(rc),x->One(corps)*Size(Intersection(u,List(w,y->x*y))));
            l[iu][iv]:=v;
            l[iv][iu]:=v;
        od;
    fi;
od;
Print("\n");
za:=Algebra(corps,l);
Print("Dimension of ZkG \t= ",Dimension(za),"\n");
radza:=RadicalOfAlgebra(za);
Print("Dimension of JZkG \t= ",Dimension(radza),"\n");
bradza:=Basis(radza);
vbradza:=BasisVectors(bradza);
vbr:=vbradza;
#
# Computing the powers of the radical of the center
i:=1;
repeat
    i:=i+1;
    l:=Set(List(Cartesian(vbradza,vbr),x->x[1]*x[2]));
    r:=Ideal(za,l);
    br:=Basis(r);
    vbr:=BasisVectors(br);
    d:=Dimension(r);
    Print("Dimension of (JZkG)^",i,"\t= ",d,"\n");
until  d=0;
\end{verbatim}
\end{footnotesize}

\section{The centre of the mod $p$ group ring of the normaliser of the Sylow subgroup of $PSU(3,p^r)$}
\label{Nsection}

\def\sur{\overline}
\newcommand{\unitmat}[3]{\left(\begin{array}{ccc}#1&#2&#3\\0&\sur{#1}/#1&-\sur{#2}/#1\\0&0&1/\sur{#1}\end{array}\right)}
\newcommand{\dunitmat}[1]{\left(\begin{array}{ccc}#1&0&0\\0&\sur{#1}/#1&0\\0&0&1/\sur{#1}\end{array}\right)}
\newcommand{\uunitmat}[2]{\left(\begin{array}{ccc}1&#1&#2\\0&1&-\sur{#1}\\0&0&1\end{array}\right)}
\def\espace{\rule{0ex}{2ex}}
\def\espaceb{\rule{0ex}{1.5ex}}
\def\espacec{\rule{0ex}{1ex}}
\newtheorem{Nota1}[Lemma1]{{Notation}}
\newenvironment{Notation}{\begin{Nota1}\em}{\end{Nota1}}
\def\dsp{\displaystyle}
\def\rmL{{\rm L}}
\def\Id{{\rm Id}}
Recall that we denote by $S$ a Sylow $p$-subgroup of the projective special unitary group $G=PSU(3,q)$ over the field with $q^2$ elements, where $q=p^r$, and by $N$ the normaliser of $S$ in $G$. In this section, we determine the ring structure of the center $ZkN$ of the group algebra $kN$, where $k$ is any commutative ring.

\begin{Notation} If $x\in N$, we denote by $x^+\in ZkN$ the sum of the conjugates of~$x$ in $N$.
\end{Notation}

Then the elements $x^+$, for $x$ in a set of representatives of conjugacy classes of~$N$, form a $k$-basis of $ZkN$.\par
Let $V$ be a three dimensional vector space over the field $\F_{q^2}$, with basis $B$. We endow $V$ with a non degenerate hermitian product, and without loss of generality, we assume that the matrix of this product in $B$ is equal to
$$\left(\begin{array}{ccc}0&0&1\\0&1&0\\1&0&0\end{array}\right)\;\;\;.$$

\begin{Notation} For $x\in \F_{q^2}$, we set $\sur{x}=x^q$. Then the map $x\mapsto\sur{x}$ is the automorphism of order 2 of the extension $\F_{q^2}/\F_q$. We also set
$$\Psi=\{x\in\F_{q^2}^\times\mid x\sur{x}=1\}$$
Let $\omega$ be a non zero element of $\F_{q^2}$ such that $\omega+\sur{\omega}=0$, and $\tau$ be an element of $\F_{q^2}$ such that $\tau+\sur{\tau}=-1$.
\end{Notation}

It follows from \cite[II Satz 10.12, page 242]{Huppert} that we can suppose that the group $N$ is equal to the image in $G$ of the group of matrices of the form
$$M(a,b,c)=\unitmat{a}{b}{c}\;\;\;,$$
where $(a, b, c)$ belongs to the set
$$\mathcal{Q}=\{(a,b,c)\in \F_{q^2}^\times\times (\F_{q^2})^2\mid b\sur{b}+a\sur{c}+c\sur{a}=0\}\;\;\;.$$

\begin{Lemma}\label{MversushatM}
For $(a,b,c)\in\mathcal{Q}$, let $\hat{M}(a,b,c)$ denote the image of $M(a,b,c)$ in $N$. Then if $(a',b',c')\in\mathcal{Q}$, we have that $\hat{M}(a,b,c)=\hat{M}(a',b',c')$ if and only if there is $\lambda\in\F_{q^2}$ with $\lambda^{q-2}=1$ and $(a',b',c')=\lambda\cdot (a,b,c)$.
\end{Lemma}

\begin{proof}
$\hat{M}(a,b,c)=\hat{M}(a',b',c')$ if and only if there exists a scalar $\lambda\in \F_{q^2}$ such that
$$\unitmat{a'}{b'}{c'}=\lambda\unitmat{a}{b}{c}\;\;\;.$$
Equivalently $(a',b',c')=\lambda(a,b,c)$ and $\sur{\lambda}/\lambda=\lambda$, i.e. $\lambda^{q-2}=1$.
\end{proof}

For two non zero integers $s,t$ denote by $(s,t)$ their greatest common divisor.
Observe that $(q-2,q^2-1)=(3,q+1)$, to motivate the following:

\begin{Notation} We set $\gamma=(3,q+1)$, and put $$\Gamma=\{\lambda\in\F_{q^2}\mid\lambda^\gamma=1\}\leq\Psi\mbox{ as well as }\rmL=\{a^\gamma\mid a\in\F_{q^2}^\times\}\;\;\;.$$
\end{Notation}

With this notation, the group $N$ has order $q^3(q^2-1)/\gamma$. It is equal to the semidirect product of the group $S$, consisting of the elements $\hat{M}(1,b,c)=\uunitmat{b}{c}$, where $b$ and $c$ are elements of $\F_{q^2}$ such that $b\sur{b}+c+\sur{c}=0$, by the cyclic group $C$ of order $(q^2-1)/\gamma$ consisting of the elements $\hat{M}(a,0,0)=\dunitmat{a}$, for $a\in \F_{q^2}^\times/\Gamma$.
%where $\Gamma=\{\lambda\in\F_{q^2}\mid\lambda^\gamma=1\}$.\par

\begin{Lemma}\label{calculs} \begin{enumerate}
\item Let $(a,b,c)$ and $(x,y,z)$ be elements of $\mathcal{Q}$. Then
$$M(a,b,c)M(x,y,z)=M\Big(ax,ay+\fracb{b\sur{x}}{x},az-\fracb{b\sur{y}}{x}+\frac{c}{\sur{x}}\Big)\;\;\;.$$
\item Let $(a,b,c)\in \mathcal{Q}$. Then $M(a,b,c)^{-1}=M\Big(\fracb{1_{\espaceb}}{\espaceb a},-\fracb{b_{\espaceb}}{\espace \sur{a}},\bar{c}\Big)$.
\item Let $(a,b,c)$ and $(x,y,z)$ be elements of $\mathcal{Q}$. Then
$$M(a,b,c)M(x,y,z)M(a,b,c)^{-1}=M\left(x,\fracb{ab}{\sur{a}}\Big(\fracb{\sur{x}}{x}-x\Big)+\fracb{a^2}{\sur{a}}y,t\right)\;\;\;,$$
where $t=a\sur{c}x+\fracb{\sur{a}c_{\espaceb}}{\espace \sur{x}}+ay\sur{b}-\fracb{\sur{a}b\sur{y}_{\espaceb}}{\espace x}+\fracb{b\sur{b}\sur{x}_{\espaceb}}{\espace x}+a\sur{a}z$.
\end{enumerate}
\end{Lemma}

\begin{proof} All the assertions follow from straightforward computations.
\end{proof}

\begin{Prop}\label{classes}\begin{enumerate}
\item Let $(x,y,z)$ and $(x',y',z')$ be elements of $\mathcal{Q}$. If $\hat{M}(x,y,z)$ and $\hat{M}(x',y',z')$ are conjugate in $N$, then $x^{-1}x'\in \Gamma$.
\item The elements $\hat{M}(x,0,0)$, for $x\in \F_{q^2}/\Gamma$, lie in distinct conjugacy classes of $N$.
\item Let $(x,y,z)\in\mathcal{Q}$. Then if $x\notin \Gamma$, the element $\hat{M}(x,y,z)$ of $N$ is conjugate to an element of the form $\hat{M}(x,0,xu\omega)$, for some $u\in \F_{q}$.
\item Let $x\in \F_{q^2}^\times$ and $u\in \F_{q}$. Then if $x\sur{x}\neq 1$, the element $\hat{M}(x,0,xu\omega)$ of $N$ is conjugate to $\hat{M}(x,0,0)$. If $x\sur{x}=1$, and if $u\neq 0$, then the element $\hat{M}(x,0,xu\omega)$ is conjugate to $\hat{M}(x,0,x\omega)$, and not conjugate to $\hat{M}(x,0,0)$.
\item If $(1,y,z)\in\mathcal{Q}$, then either $y\neq 0$ and there exists $u\in\F_q$ such that $z=y\sur{y}(\tau+u\omega)$, or $y=0$ and there exists $u\in\F_q$ such that $z=u\omega$. Moreover, if $(1,y',z')\in\mathcal{Q}$ and if $\hat{M}(1,y',z')$ and $\hat{M}(1,y,z)$ are conjugate in $N$, then $y$ and $y'$ are both non zero, or both equal to 0.
\item If $(1,y,z)$ and $(1,y',z')$ are in $\mathcal{Q}$, and if $y$ and $y'$ are both non zero, then $\hat{M}(1,y',z')$ and $\hat{M}(1,y,z)$ are conjugate in $N$ if and only if $y'^{(q^2-1)/\gamma}=y^{(q^2-1)/\gamma}$ in $\F_{q^2}^\times$, i.e. if $y'/y\in \rmL$. In particular $M(1,y,z)$ is conjugate to $M(1,y,y\sur{y}\tau)$.
%\item For $u,u'\in \F_q$, the elements $\hat{M}(1,0,u\omega)$ and $\hat{M}(1,0,u'\omega)$ are conjugate in $N$ if and only if $u$ and $u'$ are both non zero or both equal to 0. If $u\neq 0$, then $\hat{M}(1,0,u\omega)$ is conjugate to $\hat{M}(1,0,\omega)$, and if $u=0$, then $\hat{M}(1,0,u\omega)=\hat{M}(1,0,0)$ is the identity element of $N$.
\end{enumerate}
\end{Prop}

\begin{proof} Assertion (1) follows from Assertion (3) of Lemma~\ref{calculs}: if
$$\hat{M}(x',y',z')=\hat{M}\left(x,\fracb{ab}{\sur{a}}\Big(\fracb{\sur{x}}{x}-x\Big)+\fracb{a^2}{\sur{a}}y,t\right)\;\;\;,$$
then there exists $\lambda\in\Gamma$ such that $x'=\lambda x$ by Lemma~\ref{MversushatM}.

Assertion (2) is a straightforward consequence of Assertion (1).

For Assertion (3), we use Assertion (3) of Lemma~\ref{calculs} again: since $x\notin\Gamma$, we have $\fracb{\sur{x}}{x}\neq x$, and we can set $a=1$, $b=-\frac{\dsp y_{\espaceb}}{\frac{\dsp\sur{x}^{\espacec}}{\dsp x\espaceb}-{\dsp x}}$, and $c=b\sur{b}\tau$. Then $(a,b,c)\in \mathcal{Q}$ and $M(a,b,c)M(x,y,z)M(a,b,c)^{-1}$ is of the form $M(x,0,t)$, for some $t\in \F_{q^2}$. In particular $(x,0,t)\in\mathcal{Q}$, hence $x\sur{t}+t\sur{x}=0$. In other words $t=vx$ with $v+\sur{v}=0$. Then $v=u\omega$ and $u=\sur{u}$, that is $u\in\F_q$.

For Assertion (4), we have to decide when two elements of the form $n=\hat{M}(x,0,xu\omega)$ and $n'=\hat{M}(x',0,x'u'\omega)$ are conjugate in $N$, where $x, x'\notin\Gamma$, and $u,u'\in\F_q$. By Assertion~(1), we can assume that $x=x'$, and then $n$ and $n'$ are conjugate if and only if there exists $(a,b,c)\in\mathcal{Q}$ such that
$$M(a,b,c)M(x,0,xu\omega)M(a,b,c)^{-1}=M(x,0,xu'\omega)\;\;\;.$$
By Assertion (3) of Lemma~\ref{calculs}, we have $\fracb{ab}{\sur{a}}\Big(\fracb{\sur{x}}{x}-x\Big)=0$, hence $b=0$. Now $xu'\omega$ is equal to the element $t$ of Lemma~\ref{calculs}, in the case $y=b=0$ and $z=xu\omega$, that is
$$xu'\omega=a\sur{c}x+\fracb{\sur{a}c_{\espaceb}}{\espace \sur{x}}+a\sur{a}xu\omega\;\;\;.$$
Moreover $a\sur{c}+c\sur{a}=0$, since $(a,0,c)\in\mathcal{Q}$. So there exists $v\in\F_q$ such that $c=av\omega$. This gives
$$xu'\omega=-a\sur{a}xv\omega+\fracb{a\sur{a}xv\omega_{\espaceb}}{\espace \sur{x}}+a\sur{a}xu\omega\;\;\;,$$
or equivalently
$$u'=a\sur{a}\Big(u-v(1-\fracb{1}{x\sur{x}})\Big)\;\;\;.$$
Thus $n$ and $n'$ are conjugate in $N$ if and only if there exist $a\in\F_{q^2}^\times$ and $v\in\F_q$ such that $u'=a\sur{a}\Big(u-v(1-\fracb{1}{x\sur{x}})\Big)$.
If $x\sur{x}\neq 1$, then we can take $a=1$ and $v=\frac{\dsp u-u'}{{\dsp 1 -}\frac{\dsp 1\espace}{\dsp x\sur{x}\espace}}$, so $n$ and $n'$ are conjugate. And if $x\sur{x}=1$, then $n$ and $n'$ are conjugate if and only if there exists $a\in\F_{q^2}^\times$ such that $u'=a\sur{a}u$, or equivalently, if there exists $\lambda\in\F_q^\times$ such that $u'=\lambda u$. So either $u=u'=0$, or $u$ and $u'$ are both non zero. This completes the proof of Assertion (4).

For Assertion (5), assume that $(1,y,z)\in\mathcal{Q}$. Then $y\sur{y}+z+\sur{z}=0$. If $y\neq 0$, set $v=\fracb{z}{y\sur{y}\espace}-\tau$. Then $v+\sur{v}=0$, so there exists $u\in\F_q$ such that $v=u\omega$, thus $u=y\sur{y}(\tau+u\omega)$. And if $y=0$, then $z+\sur{z}=0$, so $z=u\omega$ for some $u\in \F_q$.

Now by Assertion (3) of Lemma~\ref{calculs}, for $(1,y,z)$ and $(1,y',z')$ in $\mathcal{Q}$, the elements $n=\hat{M}(1,y,z)$ and $n'=\hat{M}(1,y',z')$ are conjugate in $N$ if and only if there exists $(a,b,c)\in\mathcal{Q}$ such that
$$y'=\fracb{a^2}{\sur{a}}y\;\;\hbox{and}\;\;z'=a\sur{c}+\sur{a}c+ay\sur{b}-\sur{a}b\sur{y}+b\sur{b}+a\sur{a}z\;\;\;,$$
that is
$$y'=\fracb{a^2}{\sur{a}}y\;\;\hbox{and}\;\;z'=ay\sur{b}-\sur{a}b\sur{y}+a\sur{a}z\;\;\;,$$
In particular $y$ is non zero if and only if $y'$ is non zero. Assertion (5) follows.

Assume now that both $y$ and $y'$ are non zero. If $n$ and $n'$ are conjugate, then there exists $a\in\F_{q^2}^\times$ such that $y'=\fracb{a^2}{\sur{a}}y=a^{2-q}y$. It follows that $y'/y$ belongs to the subgroup of $\F_{q^2}^\times$ consisting of $(q-2)$-th powers, i.e. the subgroup of $\gamma$-th powers, i.e. the unique subgroup of order $(q^2-1)/\gamma$ of $\F_{q^2}^\times$. Equivalently $(y'/y)^{(q^2-1)/\gamma}=1$. Conversely, suppose that there exists $a\in\F_{q^2}^\times$ such that $y'=\fracb{a^2}{\sur{a}}y$. There are elements $u$ and $u'$ of $\F_q$ such that $z=y\sur{y}(\tau + u\omega)$ and $z'=y'\sur{y'}(\tau + u'\omega)$. If we can find $b$ and $c$ such that $(a,b,c)\in\mathcal{Q}$ and $z'=a\sur{b}y-\sur{a}b\sur{y}+a\sur{a}z$, then $n$ and $n'$ are conjugate in $N$. This can also be written
$$a\sur{a}y\sur{y}(\tau+u'\omega)=a\sur{b}y-\sur{a}b\sur{y}+a\sur{a}y\sur{y}(\tau +u\omega)\;\;\;,$$
or equivalently
$$\fracb{1}{\omega}\Big(\fracb{\sur{b}}{\sur{a}\sur{y}}-\fracb{b}{ay}\Big)=u'-u\;\;\;.\leqno{(*)}$$
Now the map $b\mapsto \fracb{1}{\omega}\Big(\fracb{\sur{b}}{\sur{a}\sur{y}}-\fracb{b}{ay}\Big)$ is a non zero $\F_q$-linear map from $\F_{q^2}$ to $F_{q}$. Hence it is surjective, and there exists $b\in\F_{q^2}$ such that (*) holds. Now we set $c=\fracb{b\sur{b}}{\sur{a}}\tau$, and then $(a,b,c)\in\mathcal{Q}$, and the elements $n$ and $n'$ are conjugate in $N$. This proves Assertion (6), and completes the proof of Proposition~\ref{classes}.
\end{proof}

\begin{Cor} The set
$$E=\big\{\hat{M}(x,0,0)\mid x\in\F_{q^2}^\times/\Gamma\big\}\bigsqcup\big\{\hat{M}(x,0,x\omega)\mid x\in\Psi/\Gamma\big\}\bigsqcup\big\{\hat{M}(1,y,y\sur{y}\tau)\mid y\in\F_{q^2}^\times/\rmL\big\}$$
is a set of representatives of conjugacy classes of $N$. In particular, there are $\fracb{{q^2+q}_{\espaceb}}{\gamma\espaceb}+\gamma$ conjugacy classes in $N$.\end{Cor}

\begin{proof} Indeed, by Proposition~\ref{classes}, the set $E$ is a set of representatives of conjugacy classes of~$N$. Its cardinality is
$$|E|=\fracb{q^2-1}{\gamma}+\fracb{q+1}{\gamma}+\gamma=\fracb{{q^2+q}_{\espaceb}}{\gamma\espaceb}+\gamma\;\;\;.$$
\end{proof}

\begin{Notation}\
\begin{itemize}
\item For $x\in \F_{q^2}^\times$, we set $d_x=\hat{M}(x,0,0)$ and $D_x=d_x^+\in ZkN$.
\item For $x\in \Psi$, we set $t_x=\hat{M}(x,0,x\omega)$ and $T_x=t_x^+$.
\item For $y^{\espace}\in \F_{q^2}^\times$, we set $u_y=\hat{M}(1,y,y\sur{y}\tau)$ and $U_y=u_y^+$.
\end{itemize}
\end{Notation}

\begin{Prop}\label{tailles}\
\begin{enumerate}
\item For $x\in \F_{q^2}^\times-\Psi$,
$$d_x^N=\big\{\hat{M}(x,y,z)\mid y,z\in\F_{q^2},\;y\sur{y}+x\sur{z}+z\sur{x}=0\big\}\;\;\;.$$
In particular $|d_x^N|=q^3$.
\item For $x\in \Psi-\Gamma$,
$$d_x^N=\big\{\hat{M}\big(x,y(\sur{x}^2-x),y\sur{y}(\sur{x}^2-x)\big)\mid y\in \F_{q^2}\big\}\;\;\;.$$
In particular $|d_x^N|=q^2$.
\item For $x\in \Gamma$, the element $d_x$ is the identity element of $N$, and $|d_x^N|=1$.
\item For $x\in\Psi$, the conjugacy class of $t_x$ in $N$ has cardinality $q^2(q-1)$ if $x\notin\Gamma$, and $q-1$ otherwise. The conjugacy class of $T_1$ consists of the elements $\hat{M}(1,0,\lambda\omega)$, for $\lambda\in\F_{q}^\times$.
\item For $x\in \F_{q^2}^\times$,
$$u_x^N=\{\hat{M}(1,v,v\sur{v}\tau+\lambda\omega)\mid v\in x\rmL,\;\lambda\in\F_q\}\;\;\;.$$
In particular $|u_x^N|=\fracb{q(q^2-1)}{\gamma\espaceb}$.
\end{enumerate}
\end{Prop}

\begin{proof} It follows from Proposition~\ref{classes} that if $(x,y,z)\in\mathcal{Q}$ and $x\sur{x}\neq 1$, then $\hat{M}(x,y,z)$ is conjugate to $d_x$, and that conversely, any conjugate of $d_x$ in $N$ is of the form $\hat{M}(x,y,z)$, for some elements $y,z\in\F_{q^2}$ such that $(x,y,z)\in\mathcal{Q}$. This proves Assertion (1).\par
Now let $(a,b,c)$ and $(x,y,z)$ be elements of $\mathcal{Q}$. By Assertion (1) of Lemma~\ref{calculs}, comparing the diagonal elements in the product in the two possible orders, the elements $\hat{M}(a,b,c)$ and $\hat{M}(x,y,z)$ commute if and only if
$$ay+\fracb{b\sur{x}}{x}=xb+\fracb{y\sur{a}}{a}\;\;\hbox{and}\;\;az-\fracb{b\sur{y}}{x}+\frac{c}{\sur{x}}=xc-\fracb{y\sur{b}}{a}+\frac{z}{\sur{a}}\;\;\;\leqno{(**)}$$
\begin{itemize}
\item If $y=z=0$, this gives $\fracb{b\sur{x}}{x}=xb$ and $\fracb{c}{\sur{x}}=xc$. If moreover $x\sur{x}=1$ but $x^2\neq\sur{x}$, then $b=0$, but $a$ and $c$ are arbitrary, only subject to $a\sur{x}+c\sur{a}=0$. In this case the centraliser of $d_x$ in $N$ has cardinality $\fracb{q(q^2-1)}{\gamma}$, and the conjugacy class of $d_x$ in $N$ has cardinality $q^2$. Now Assertion (2) follows from the fact that the elements
$$\hat{M}(1,y,y\sur{y}\tau)\hat{M}(x,0,0)\hat{M}(1,y,y\sur{y}\tau)^{-1}=\hat{M}\big(x,y(\sur{x}^2-x),y\sur{y}(\sur{x}^2-x)\big)\;\;\;,$$
for $y\in \F_{q^2}$, are all distinct.

Finally if $x^2=\sur{x}$, then $x\in\Gamma$, so $d_x$ is the identity element of $N$, and Assertion~(3) follows.

\item If $x\in\Psi$, $y=0$ and $z=x\omega$, then the relations (**) give
$$b\sur{x}^2=xb\;\;\hbox{and}\;\;ax\omega+cx=xc+\frac{x\omega}{\sur{a}}\;\;\;,$$
that is $b(x-\sur{x}^2)=0$ and $a\sur{a}=1$. If $x\neq\sur{x}^2$, i.e. if $x\notin\Gamma$, this is equivalent to $b=0$ and $a\sur{a}=1$. Then $c$ is arbitrary, only subject to $a\sur{x}+c\sur{a}=0$. In this case the centraliser of $t_x$ in $N$ has cardinality $\fracb{q(q+1)}{\gamma}$, and the conjugacy class of $t_x$ in $N$ has cardinality $q^2(q-1)$. Now if $x^2=\sur{x}$, the only condition left is $a\sur{a}=1$, so the centraliser of $t_x$ in $N$ has cardinality $q^3$ (it is equal to $S$), and the conjugacy class of $t_x$ in $N$ has cardinality $q-1$. Moreover, by Lemma~\ref{calculs}, the conjugates of $t_1=\hat{M}(1,0,\omega)$ are the elements $\hat{M}(1,0,a\sur{a}\omega)$, for $a\in\F_{q^2}^\times$. This completes the proof of Assertion (4).

\item If $x=1$, $y\in\F_{q^2}^\times$, and $z=y\sur{y}\omega$, then the relations (**) give
$$ay=\fracb{y\sur{a}}{a}\;\;\hbox{and}\;\;ay\sur{y}\omega-b\sur{y}=-\fracb{y\sur{b}}{a}+\frac{y\sur{y}\omega}{\sur{a}}\;\;\;.$$
Since $y\neq 0$, the first relation gives $a^2=\sur{a}$, i.e. $a\in\Gamma$, so we can assume $a=1$ by Lemma~\ref{MversushatM}. Now the second relation reads $b\sur{y}=y\sur{b}$, i.e. $b=uy$, for $u\in\F_q$. Since $c$ is subject to $c+\sur{c}+b\sur{b}=0$, it follows that the centraliser of $u_y$ in $N$ has cardinality $q^2$, and the conjugacy class of $u_y$ in $N$ has cardinality $\fracb{q(q^2-1)}{\gamma\espaceb}$.\par
Now Assertion (5) follows from the fact that by Proposition~\ref{classes}, the element $\hat{M}(1,v,v\sur{v}\tau+\lambda\omega)$, for $v\in x\rmL$ and $\lambda\in \F_q$, is conjugate to $u_x$, and that there are $\fracb{q(q^2-1)}{\gamma\espaceb}$ such elements in $N$.
\end{itemize}
\end{proof}

We recall the following well known fact (cf e.g.~\cite[(9.28)]{CR}):

\begin{Lemma}\label{classique} Let $G$ be a finite group and $k$ be a commutative ring. For $x\in G$, let $x^+\in ZkG$ denote the sum of the elements of the conjugacy class $x^G$ of $x$ in $G$. Then for $x,y\in G$
$$x^+\cdot y^+=\sum_{z\in[G]}m_{x,y}^zz^+\;\;\;,$$
where $[G]$ denotes a set of representatives of conjugacy classes of $G$, and
$$m_{x,y}^z=|\{(x',y')\in x^G\times y^G\mid x'y'=z\}|\;\;\;.$$
Clearly $m_{x,y}^z=m_{y,x}^z$ and $m_{x^{-1},y^{-1}}^{z^{-1}}=m_{x,y}^z$ for any $x,y,z\in G$, but since
$$m_{x,y}^z|z^G|=|\{(x',y',z')\in x^G\times y^G\times z^G\mid x'y'=z'\}|\;\;\;,$$
we have also $m_{x,y}^z|z^G|=m_{z,y^{-1}}^x|x^G|=m_{z,x^{-1}}^y|y^G|$.
\end{Lemma}

Observe that $Z(kN)=k\otimes_\Z Z(\Z N)$ and hence we may and will suppose for the rest of this section that $k=\Z$, unless otherwise stated.

\begin{Prop}\label{premiers produits} \begin{enumerate}
\item Let $x\in F_{q^2}^\times-\Psi$ and $y\in\F_{q^2}^\times$ such that $xy\notin \Psi$. Then
$$D_xD_y=\left\{\begin{array}{cl}q^3D_{xy}&\hbox{if}\;y\notin \Psi\\q^2D_{xy}&\hbox{if}\;y\in \Psi\;\;\;.\end{array}\right.$$
\item Let $x\in F_{q^2}^\times-\Psi$ and $y\in\Psi$. Then
$$D_xT_y=\left\{\begin{array}{cl}q^2(q-1)D_{xy}&\hbox{if}\;y\notin\Gamma\\(q-1)D_{xy}&\hbox{if}\; y\in\Gamma\;\;\;.\end{array}\right.$$
\item Let $x\in F_{q^2}^\times-\Psi$ and $y\in\F_{q^2}^\times$. Then $D_xU_y=\fracb{q(q^2-1)}{\gamma\espaceb}D_{x}$.
\end{enumerate}
\end{Prop}

\begin{proof} The three assertions follow from the fact that  the product of an element in the conjugacy class of $r=\hat{M}(x_1,y_1,z_1)$ of $N$ and an element in the conjugacy class of $s=\hat{M}(x_2,y_2,z_2)$ of $N$ is an element of the form $\hat{M}(x_1x_2,\alpha,\beta)$, for some $\alpha$ and $\beta$ in $\F_{q^2}$. In each assertion, the assumption implies that all these elements are in the conjugacy class of $t=d_{x_1x_2}$, since $x_1x_2\in\F_{q^2}-\Psi$. It follows that there exists an integer $m$ such that $r^+s^+=mD_{x_1x_2}$.\par
Now the augmentation map $\varepsilon:kN\to k$ restricts to a ring homomorphism $ZkN\to k$, sending $x^+$ to $|x^G|$. Hence $|r^N||s^N|=m|t^N|$. For the three assertions, we can assume that $r=d_x$ and $x\notin\Psi$, thus $|r^N|=q^3$. Similarly $t=d_{xy}$ for Assertions (1) and (2), and $xy\notin\Psi$, so $|t^N|=q^3$. For Assertion (3), we have $t=d_x$, so $|t^N|=q^3$ again. It follows that the integer $m$ is equal to $|s^N|$, and $s=d_y$ in Assertion (1), $s=t_y$ in Assertion (2), and $s=u_y$ in Assertion (3). Now Proposition~\ref{premiers produits} follows from the values of the cardinalities $|s^N|$ given by Proposition~\ref{tailles}.
\end{proof}

\begin{Prop}\label{DxDyprodinPsi} Let $x,y\in\F_{q^2}-\Psi$, such that $xy\in\Psi-\Gamma$. Then $D_xD_y=q^3D_{xy}+q^3T_{xy}$.
\end{Prop}

\begin{proof}
Any element in the product $d_x^N\cdot d_y^N$ is of the form $\hat{M}(xy,\alpha,\beta)$, for some $\alpha,\beta\in\F_{q^2}$. It follows that there are integers $a$ and $b$ such that $D_xD_y=aD_{xy}+bT_{xy}$. Setting $z=xy$, the integer $a$ is equal to $m_{d_x,d_y}^{d_z}$. Thus $a|d_z^N|=m_{d_z,d_y^{-1}}^{d_x}|d_x^N|$, by Lemma~\ref{classique}. But by Proposition~\ref{premiers produits}, we have $D_zD_{y^{-1}}=q^2D_{x}$, so $m_{d_z,d_y^{-1}}^{d_x}=q^2$. It follows that $a|z^N|=aq^2=q^2q^3$, thus $a=q^3$. Taking augmentation gives
$$\varepsilon(D_xD_y)=q^6=a\varepsilon(D_z)+b\varepsilon(T_z)=aq^2+bq^2(q-1)\;\;\;.$$
It follows that $b=\fracb{{q^6-q^5}_{\espaceb}}{q^{\espaceb 2}(q-1)}=q^3$.
\end{proof}

\begin{Prop} \label{DxT1} Let $x\in\Psi$. Then $D_xT_1=T_x$.
\end{Prop}

\begin{proof}
If $x\in\Gamma$, there is nothing to prove, because $D_x$ is equal to the identity, in this case. If $x\notin\Gamma$, then $D_xT_1$ is a sum of elements of the form $\hat{M}(x,\alpha,\beta)$, so there are natural integers $a$ and $b$ such that $D_xT_1=aD_x+bT_x$. Taking augmentation of this equality gives $q^2(q-1)=aq^2+bq^2(q-1)$, that is $q-1=a+b(q-1)$. Since the product $d_xt_1$ is equal to~$t_x$, it follows that $b>0$. Hence $b=1$ and $a=0$.
\end{proof}

\begin{Prop} \label{DxUy} Let $x\in\Psi-\Gamma$, and $y\in \F_{q^2}^\times$. Then $D_xU_y=\fracb{q^2-1}{\gamma\espaceb}(D_x+T_x)$.
\end{Prop}

\begin{proof} Again $D_xU_y$ is a sum of elements of $N$ of the form $\hat{M}(x,\alpha,\beta)$. Hence there are natural integers $a$ and $b$ such that $D_xU_y=aD_x+bT_x$. The integer $a$ is equal to $m_{d_x,u_y}^{d_x}$, i.e.
$$a=|\{(d',u')\in d_x^N\times u_y^N\mid d'u'=d_x\}|\;\;\;.$$
By Proposition~\ref{tailles}, the element $d'\in d_x^N$ is equal to $\hat{M}\big(x,w(\sur{x}^2-x),w\sur{w}(\sur{x}^2-x)\big)$, for $w\in\F_{q^2}$, and the element $u'$ is equal to $\hat{M}(1,v,v\sur{v}\tau +\lambda\omega)$, for $v\in x\rmL$ and $\lambda\in\F_q$. Now
$$d'u'=\hat{M}\big(x,xv+w(\sur{x}^2-x),x(v\sur{v}\tau+\lambda\omega)-w\sur{v}(\sur{x}^2-x)+w\sur{w}(\sur{x}^2-x)\big)\;\;\;.$$
This is equal to $d_x$ if and only if
$$xv+w(\sur{x}^2-x)=0\;\;\hbox{and}\;\;x(v\sur{v}\tau+\lambda\omega)-w\sur{v}(\sur{x}^2-x)+w\sur{w}(\sur{x}^2-x)=0\;\;\;.$$
Since $x\notin\Gamma$, the first relation gives $w=\fracb{v}{1-\sur{x}^3}$. Multiplying by $\sur{x}$, the second one reads
$$v\sur{v}\tau+\lambda\omega -w\sur{v}(\sur{x}^3-1)+w\sur{w}(\sur{x}^3-1)=0\;\;\;.$$
This gives
$$v\sur{v}\tau+\lambda\omega+v\sur{v}-\fracb{v\sur{v}}{1-x^3}=0\;\;\;,$$
that is
$$\lambda=\fracb{1}{\omega}\big(\sur{\tau}+\fracb{1}{1-x^3}\big)\;\;\;.$$
This defines an element $\lambda$ of $\F_q$, since $\tau+\sur{\tau}=-1$ and
$$\fracb{1}{1-x^3}+\fracb{1}{1-\sur{x}^3}=\fracb{2-x^3-\sur{x}^3}{(1-x^3)(1-\sur{x}^3)}=1\;\;\;.$$
In other words $w$ and $\lambda$ are determined by $v\in x\rmL$, which may be chosen  arbitrarily.  It follows that $a=\fracb{q^2-1}{\gamma\espaceb}$.\par
Now applying the augmentation to the relation $D_xU_y=aD_x+bT_x$ gives
$$q^2\fracb{q(q^2-1)}{\gamma\espaceb}=aq^2+bq^2(q-1)\;\;\;.$$
It follows that
$$\fracb{q(q^2-1)}{\gamma\espaceb}=\fracb{q^2-1}{\gamma\espaceb}+b(q-1)\;\;\;,$$
hence $b=\fracb{q^2-1}{\gamma\espaceb}$.
\end{proof}
\pagebreak[3]
\begin{Prop} \label{DxDxinv}\
\begin{enumerate}
\item Let $x\in\Psi-\Gamma$. Then
$$D_xD_{x^{-1}}=q^2\Id+q\sum_{y\in\F_{q^2}^\times/\rmL}U_y\;\;\;.$$
\item Let $x\in\F_{q^2}^\times-\Psi$. Then
$$D_xD_{x^{-1}}=q^3\Id+q^3T_1+q^3\sum_{y\in\F_{q^2}^\times/\rmL}U_y\;\;\;.$$
\end{enumerate}
\end{Prop}

\begin{proof} For $x\in\F_{q^2}^\times$, the product $D_xD_{x^{-1}}$ is a sum of elements of the form $\hat{M}(1,\alpha,\beta)$ of $N$. So there are integers $a,b,c_y\in\N$, for $y\in\F_{q^2}^\times/\rmL$ such that
$$D_xD_{x^{-1}}=a\Id+bT_1+\sum_{y\in\F_{q^2}^\times/\rmL}c_yU_y\;\;\;.\leqno{(***)}$$
Then $a=m_{d_x,d_{x^{-1}}}^{\Id}=|\{(d',d'')\in d_x^N\times d_{x^{-1}}^N\mid d' d''=\Id\}|=|d_x^N|$. Thus $a=q^2$ if $x\in\Psi-\Gamma$, and $a=q^3$ if $x\in\F_{q^2}^\times-\Psi$.\par
On the other hand, by Lemma~\ref{classique}, for $y\in\F_{q^2}^\times$,
$$c_y|u_y^N|=m_{d_x,d_{x^-1}}^{u_y}|u_y^N|=m_{u_y,d_x}^{d_x}|d_x^N|\;\;\;$$
\begin{itemize}
\item If $x\in\Psi-\Gamma$, then $m_{u_y,d_x}^{d_x}=\fracb{q^2-1}{\gamma\espaceb}$, by Proposition~\ref{DxUy}. It follows that
$$c_y\fracb{q(q^2-1)}{\gamma\espaceb}=\fracb{q^2-1}{\gamma\espaceb}q^2\;\;\;,$$
hence $c_y=q$.\par
Applying augmentation to equation (***), we get $q^2q^2=a+b(q-1)+q\cdot \gamma q\fracb{q^2-1}{\gamma\espaceb}$. This gives $b(q-1)=q^4-q^2-q^2(q^2-1)=0$, which proves Assertion (1).
\item If $x\in\F_{q^2}^\times-\Psi$, then $m_{u_y,d_x}^{d_x}=\fracb{q(q^2-1)}{\gamma\espaceb}$ by Proposition~\ref{premiers produits}. Thus $c_y=q^3$ in this case. Applying augmentation to equation (***) gives
$$q^3\cdot q^3=q^3+b(q-1)+q^3\gamma\cdot \fracb{q(q^2-1)}{\gamma\espaceb}\;\;\;,$$
that is $b(q-1)=q^6-q^3-q^4(q^2-1)=q^3(q-1)$, hence $b=q^3$, which proves Assertion~(2).
\end{itemize}
\end{proof}

\begin{Prop} \label{DxDyPsi} Let $x,y\in\Psi-\Gamma$ such that $xy\notin\Gamma$. Then $D_xD_y=D_{xy}+(q+1)T_{xy}$.
\end{Prop}

\begin{proof} The product $D_xD_y$ is a sum of elements of $N$ of the form $\hat{M}(xy,\alpha,\beta)$, so there are integers $a$ and $b$ such that $D_xD_y=aD_{xy}+bT_{xy}$. The integer $a$ is the number of pairs $(d',d'')$ in $d_x^N\times d_y^N$ such that $d'd''=d_{xy}$.

By Proposition~\ref{tailles}, the class $d_x^N$ consists of the elements $\hat{M}\big(x,\alpha(\sur{x}^2-x),\alpha\sur{\alpha}(\sur{x}^2-x)\big)$, for $\alpha\in\F_{q^2}$.
Equivalently, in a form that will be more convenient for computation, it consists of the elements $d'=\hat{M}(x,u,v)$, for $u\in \F_{q^2}$ and $v=\fracb{u\sur{u}}{x^2-\sur{x}}$.
Similarly, the class $d_y^N$ consist of the elements $d''=\hat{M}(y,r,s)$,
for $r\in\F_{q^2}$ and $s=\fracb{r\sur{r}}{y^2-\sur{y}}$. Since $x\sur{x}=1=y\sur{y}$, we have
$$d'd''=\left(\begin{array}{ccc}x&u&v\\0&\sur{x}^2&-\sur{u}\,\sur{x}\\0&0&x\end{array}\right)
\left(\begin{array}{ccc}y&r&s\\0&\sur{y}^2&-\sur{r}\,\sur{y}\\0&0&y\end{array}\right)\;\;\;.$$
The product $d'd''$ is equal to $d_{xy}$ if and only if
$$xr+u\sur{y}^2=0\;\;\hbox{and}\;\;xs-u\sur{r}\,\sur{y}+vy=0\;\;\;.$$
The first equation gives $r=-u\sur{x}\,\sur{y}^2$, thus $r\sur{r}=u\sur{u}$. Now the second equation becomes
$$\fracb{xu\sur{u}}{y^2-\sur{y}}+u\sur{u}xy+\fracb{yu\sur{u}}{x^2-\sur{x}}=0\;\;\;.$$
Then either $u=0$, hence $r=s=v=0$, or
$$\fracb{x}{y^2-\sur{y}}+xy+\fracb{y}{x^2-\sur{x}}=0\;\;\;.$$
Equivalently $(x^3-1)+(x^3-1)(y^3-1)+(y^3-1)=0$, thus $x^3y^3=1$, which doesn't hold since $xy\notin\Gamma$, using the remark after Lemma~\ref{MversushatM}.

So the only pair $(d',d'')\in d_x^N\times d_y^N$ such that $d'd''=d_{xy}$ is the pair $(d_x,d_y)$. It follows that $a=1$.

Applying augmentation to the equality $D_xD_y=aD_{xy}+bT_{xy}$ now gives $q^4=q^2+bq^2(q-1)$, hence $b=q+1$
\end{proof}

\begin{Prop}\label{TxDyPsi} Let $x\in\Psi-\Gamma$ and $y\in\F_{q^2}^\times$ with $xy\notin\Gamma$. Then
$$D_xT_y=(q^2-1)D_{xy}+(q^2-q-1)T_{xy}\;\;\;.$$
\end{Prop}

\begin{proof} The product $D_xT_y$ is a sum of elements of the form $\hat{M}(xy,\alpha,\beta)$, so there are integers $a$ and $b$ such that $D_xT_y=aD_{xy}+bT_{xy}$. By Lemma~\ref{classique}, Proposition~\ref{tailles}, and Proposition~\ref{DxDyPsi}, we have
$$aq^2=m_{d_x,t_y}^{d_{xy}}|d_{xy}^N|=m_{d_{xy},d_{x}^{-1}}^{t_y}q^2(q-1)=q^2(q^2-1)\;\;\;,$$
hence $a=q^2-1$. Taking augmentation gives
$$\varepsilon(D_xT_y)=q^2q^2(q-1)=a\varepsilon(D_{xy})+b\varepsilon(T_{xy})=(q^2-1)q^2+bq^2(q-1)\;\;\;,$$
hence $b=q^2-q-1$.
\end{proof}

\begin{Prop} \label{TxT1}\begin{enumerate}
\item $T_1^2=(q-1)\Id+(q-2)T_1$.
\item If $x\in\Psi-\Gamma$, then $T_xT_1=(q-1)D_x+(q-2)T_x$.
\end{enumerate}
\end{Prop}

\begin{proof} By Proposition~\ref{tailles}, the product of any two conjugates of $t_1$ is either the identity, or again a conjugate of $t_1$. It follows that there are integers $a$ and $b$ such that $T_1^2=a\Id+bT_1$. Moreover $a$ is equal to the cardinality of the conjugacy class of $t_1$, that is $a=q-1$. Now taking augmentation gives $(q-1)^2=a+(q-1)b$, hence $b=q-2$. Now for $x\in\Psi-\Gamma$,
$$T_xT_1=D_xT_1^2=(q-1)D_x+(q-2)T_x\;\;\;,$$
since $D_xT_1=T_x$ by Proposition~\ref{DxT1}.
\end{proof}

\begin{Prop} \label{DxTxinv}
Let $x\in\Psi-\Gamma$. Then $D_xT_{x^{-1}}=q^2T_1+q(q-1)\sum_{y\in \F_{q^2}^\times/{\rm L}}\limits U_y$.
\end{Prop}

\begin{proof} Again $D_xT_{x^{-1}}$ is a sum of elements of the form $\hat{M}(1,\alpha,\beta)$, so there are integers $a,b$, and $c_y$, for $y\in \F_{q^2}^\times/{\rm L}$, such that $D_xT_{x^{-1}}=a\Id+bT_1+\sum_{y\in \F_{q^2}^\times/{\rm L}}\limits c_yU_y$. Since $t_{x^{-1}}=t_x^{-1}$, and since no conjugate of $d_x$ is a conjugate of $t_x$, we have $a=0$. Then $b=m_{d_x,t_{x^-1}}^{t_1}$, hence $b(q-1)=m_{t_1,d_x}^{t_x}q^2(q-1)=q^2(q-1)$, by Proposition~\ref{DxT1}. Hence $b=q^2$. Similarly $c_y=m_{d_x,t_{x^{-1}}}^{u_y}$, so $c_y\fracb{q(q^2-1)}{\gamma}=m_{u_y,d_{x^{-1}}}^{t_{x^{-1}}}q^2(q-1)$, hence $c_y\fracb{q(q^2-1)}{\gamma}=\fracb{q^2-1}{\gamma}q^2(q-1)$, thus $c_y=q(q-1)$.
\end{proof}

\begin{Prop} \label{TxUy} Let $x\in\Psi-\Gamma$ and $y\in\F_{q^2}$. Then $T_xU_y=\fracb{(q^2-1)(q-1)}{\gamma}(D_x+T_x)$.
\end{Prop}

\begin{proof} By Proposition~\ref{DxT1} and Proposition~\ref{DxUy}, we have that
\begin{eqnarray*}
T_xU_y=D_xT_1U_y&=&\fracb{(q^2-1)}{\gamma}(D_x+T_x)T_1\\
&=&\fracb{(q^2-1)}{\gamma}\big(T_x+(q-1)D_x+(q-2)T_x\big)\\
&=&\fracb{(q^2-1)(q-1)}{\gamma}(D_x+T_x)
\end{eqnarray*}
\end{proof}

\begin{Prop} Let $x\in\Psi-\Gamma$. Then
$$T_xT_{x^{-1}}=q^2(q-1)\Id +q^2(q-2)T_1+q(q-1)^2\sum_{y\in \F_{q^2}^\times/{\rm L}}\limits U_y\;\;\;.$$
\end{Prop}

\begin{proof} Indeed by Proposition~\ref{DxDxinv}, Proposition~\ref{DxT1}, Proposition~\ref{TxT1} and Proposition~\ref{DxTxinv}
\begin{eqnarray*}
T_xT_{x^{-1}}&=&D_xT_1D_{x^{-1}}T_1\\
&=&D_xD_{x^{-1}}T_1^2\\
&=&D_xD_{x^{-1}}\big((q-1)\Id+(q-2)T_1\big)\\
&=&D_x\big((q-1)D_{x^{-1}}+(q-2)T_{x^{-1}}\big)\\
&=&(q-1)\Big(q^2\Id+q\sum_{y\in \F_{q^2}^\times/\rmL}\limits U_y\Big)+(q-2)\Big(q^2T_1+q(q-1)\sum_{y\in \F_{q^2}^\times/{\rm L}}\limits U_y\Big)\\
&=&q^2(q-1)\Id +q^2(q-2)T_1+q(q-1)^2\sum_{y\in \F_{q^2}^\times/{\rm L}}\limits U_y\;\;\;.
\end{eqnarray*}
\end{proof}

\begin{Prop} \label{TxTy}
Let $x, y\in \Psi-\Gamma$ such that $xy\notin\Gamma$. Then
$$T_xT_y=(q-1)(q^2-q-1)D_{xy}+\big(q(q-1)^2+1\big)T_{xy}\;\;\;.$$
\end{Prop}

\begin{proof} Indeed, by Proposition~\ref{DxDyPsi}, Proposition~\ref{DxT1}  and Proposition~\ref{TxT1}
\begin{eqnarray*}
T_xT_y&=&D_xT_1D_yT_1\\
&=&\big(D_{xy}+(q+1)T_{xy}\big)\big((q-1)\Id+(q-2)T_1\big)\\
&=&(q-1)D_{xy}+(q-2)T_{xy}+(q^2-1)T_{xy}+(q-2)(q+1)\big((q-1)D_{xy}+(q-2)T_{xy}\big)\\
&=&(q-1)(q^2-q-1)D_{xy}+\big(q(q-1)^2+1\big)T_{xy}
\end{eqnarray*}
\end{proof}

\begin{Prop} \label{T1Ux} Let $x\in\F_{q^2}^\times$. Then $T_1U_x=(q-1)U_x$.
\end{Prop}

\begin{proof} The product $T_1U_x$ is a linear combination of elements of $N$ of the form $\hat{M}(1,\alpha,\beta)$. Hence there are integers $a,b$ and $c_y$, for $y\in \F_{q^2}^\times/{\rm L}$, such that
$$T_1U_x=a\Id+bT_1+\sum_{y\in \F_{q^2}^\times/{\rm L}}c_yU_y\;\;\;.\leqno(\#)$$
Observe now that $t_1$ and $u_x^{-1}$ are not conjugate in $N$, e.g. because the conjugacy class of $t_1$ has cardinality $q-1$, and the conjugacy class of $u_x$ has cardinality $\fracb{q(q^2-1)}{\gamma}\neq q-1$. It follows that $a=0$.\par
Now by Proposition~\ref{tailles}, the conjugacy class of $T_1$ consists of the elements $\hat{M}(1,0,\lambda\omega)$, for $\lambda\in\F_q^\times$, and the conjugacy class of $u_x$ consists of the elements $\hat{M}(1,v,v\sur{v}\tau+\mu\omega)$, for $v\in x{\rm L}$ and $\mu\in\F_q$. The product $\pi=\hat{M}(1,0,\lambda\omega)\hat{M}(1,v,v\sur{v}\tau+\mu\omega)$ is equal to $u_y=\hat{M}(1,y,y\sur{y}\tau)$ if and only if $v=y$ and $v\sur{v}\tau+\mu\omega+\lambda\omega=y\sur{y}\tau$. It follows that $c_y=0$ unless $y\in x{\rm L}$, i.e. unless $y{\rm L}=x{\rm L}$. If $y{\rm L}=x{\rm L}$, then $u_y$ is conjugate to $u_x$ in $N$, and we can assume that $y=x$. In this case $\pi=u_x$ if and only if $v=x$ and $\mu=-\lambda$. It follows that $c_x=q-1$.

Applying augmentation to Equation~(\#) now gives
'$$(q-1)\cdot \fracb{q(q^2-1)}{\gamma}=b(q-1)+(q-1)\cdot \fracb{q(q^2-1)}{\gamma}\;\;\;,$$
hence $b=0$.
\end{proof}

\begin{Prop} \label{UxUy}
\begin{enumerate}
\item If $3\nmid q+1$, then ${\rm L}=\F_{q^2}^\times$, and
$$U_1^2=q(q^2-1)\Id+q(q^2-1)T_1+q(q^2-2)U_1\;\;\;.$$
\item If $3\mid q+1$, then $\F_{q^2}/{\rm L}=\{{\rm L},t{\rm L},t^2{\rm L}\}$, where $t$ is any non cube element of $\F_{q^2}^\times$. Let $l=|\{v\in {\rm L}\mid 1-v\in {\rm L}\}|$, $m=|\{v\in {\rm L}\mid t-v\in {\rm L}\}|$, and $n=|\{v\in {\rm L}\mid t-v/t\in {\rm L}\}|$. Then for $x\in \F_{q^2}^\times/{\rm L}$,
\begin{eqnarray*} U_x^2&=&\fracb{q(q^2-1)}{\gamma}(\Id+T_1)+ql U_x+qm(U_{tx}+U_{t^2x})\\
U_xU_{tx}&=&qn U_{t^2x}+qm(U_x+U_{tx})\;\;\;.
\end{eqnarray*}
\end{enumerate}
\end{Prop}

\begin{proof}
By Proposition~\ref{tailles}, for $x\in\F_{q^2}^\times$, the conjugacy class of $u_x$ in $N$ consists of the elements $\hat{M}(1,v,v\sur{v}\tau+\lambda\omega)$, for $v\in x{\rm L}$ and $\lambda\in \F_q$. Since the inverse of $u_x=\hat{M}(1,x,x\sur{x}\tau)$ is $\hat{M}(1,-x,x\sur{x}\tau)$, and since $-x\in x{\rm L}$ as $-1=(-1)^\gamma\in {\rm L}$, we have that $u_x^{-1}$ is conjugate to~$u_x$.

For $x,y\in \F_{q^2}^\times$, the product $U_xU_y$ is a sum of elements of the form $\hat{M}(1,\alpha,\beta)$, hence there are integers $a, b$ and $c_{x,y}^z$, for $z\in\F_{q^2}^\times/{\rm L}$, such that
$$U_xU_y=a\Id+bT_1+\sum_{z\in\F_{q^2}^\times/{\rm L}}c_{x,y}^zU_z\;\;\;.\leqno(\#\#)$$
Note that for $x,y,z\in \F_{q^2}^\times$, we have
$$c_{x,y}^z|u_z^N|=m_{u_x,u_y}^{u_z}\fracb{q(q^2-1)}{\gamma}=
m_{u_z,u_x^{-1}}^{u_y}\fracb{q(q^2-1)}{\gamma}=
c_{z,x}^y\fracb{q(q^2-1)}{\gamma}=c_{z,x}^y|u_z^N|\;\;\;,$$
as $u_{x}^{-1}$ is conjugate to $u_x$. So $c_{x,y,z}$ is a symmetric function of $x,y,z$.

If $x{\rm L}\neq y{\rm L}$, then no conjugate of $u_x^{-1}$ is conjugate to $u_y$, so $a=0$. In this case, we also have
$$b|t_1^N|=m_{u_x,u_y}^{t_1}(q-1)=m_{t_1,u_{x}^{-1}}^{u_y}|u_y^N|\;\;\;,\leqno(\#\#\#)$$
and $m_{t_1,u_{x}^{-1}}^{u_y}=0$ by Proposition~\ref{T1Ux}. It follows that $b=0$ in this case.

If $x{\rm L}=y{\rm L}$, i.e. $U_x=U_y$, then clearly $a=|u_x^N|=\fracb{q(q^2-1)}{\gamma}$. Moreover Equation~(\#\#\#) gives $b(q-1)=(q-1)|u_y^N|$, hence $b=\fracb{q(q^2-1)}{\gamma}$.

In the case $3\nmid q+1$, we have $\gamma=1$ and ${\rm L}=\F_{q^2}^\times$. Then
$$U_1^2=q(q^2-1)(\Id+T_1)+c_{1,1}^1U_1\;\;\;.$$
Taking augmentation gives
$$\big(q(q^2-1)\big)^2=q(q^2-1)(1+q-1)+c_{1,1}^1q(q^2-1)\;\;\;,$$
hence $c_{1,1}^1=q(q^2-2)$, which completes the proof of Assertion~(1).

In the case $3\mid q+1$, then $\gamma=3$, and $\rmL$ has index 3 in $\F_{q^2}^\times$, so $\F_{q^2}^\times/{\rm L}=\{1,t{\rm L},t^2{\rm L}\}$ for any non cube element $t$ of $\F_{q^2}^\times$.

For $x,y,z\in \F_{q^2}^\times$, the product of the element $u'=\hat{M}(1,v,v\sur{v}\tau+\lambda\omega)$ in the conjugacy class of $u_x$ (where $v\in x\rmL$ and $\lambda\in \F_q$) by the element $u''=\hat{M}(1,r,r\sur{r}\tau+\mu\omega)$ in the conjugacy class of $u_y$ (where $r\in y\rmL$ and $\mu\in\F_q$) is equal to $u_z$ if and only if
$$v+r=z\;\;\hbox{and}\;\;r\sur{r}\tau+\mu\omega-v\sur{r}+v\sur{v}\tau+\lambda\omega=
z\sur{z}\tau\;\;\;.$$
The second equation determines $\mu$ once $v,r$ and $\lambda$ are known, and $\lambda$ can be chosen arbitrarily in $\F_q$, once $v$ and $r$ satisfy $v+r=z$. Hence in Equation~(\#\#), we have
$$c_{x,y}^z=q\;\big|\{v\in x{\rm L}\mid z-v\in y{\rm L}\}\big|\;\;\;.$$
In particular for any $x\in\F_{q^2}^\times$
$$c_{x,x}^x=q\;\big|\{v\in x{\rm L}\mid x-v\in x{\rm L}\}\big|=q\;\big|\{w\in {\rm L}\mid x-xw\in x{\rm L}\}\big|=l\;\;\;.$$
Similarly
$$c_{x,x}^{xt}=q\;\big|\{v\in x{\rm L}\mid xt-v\in x{\rm L}\}\big|=q\;\big|\{w\in {\rm L}\mid xt-xw\in x{\rm L}\}\big|=m\;\;\;.$$
Finally
$$c_{x,xt}^{xt^2}=q\;\big|\{v\in x{\rm L}\mid xt^2-v\in xt{\rm L}\}\big|=q\;\big|\{w\in {\rm L}\mid t^2-w\in t{\rm L}\}\big|=n\;\;\;.$$
This completes the proof, since $c_{x,y}^z$ is symmetric in $x,y,z$.
\end{proof}

\begin{Rem}
Applying augmentation to the equations of Proposition~\ref{UxUy} gives that $n=l+1$ and $n+2m=\fracb{q^2-1}{3}$. So it suffices to know $l$, and then $m$ and $n$ can be computed.

By definition $l=|\{v\in {\rm L}\mid 1-v\in {\rm L}\}|$. Since $3\mid q+1\mid q^2-1$, the field $\F_{q^2}$ contains all cubic roots of unity. Now clearly
$$l=|\{(x,y)\in \F_{q^2}^\times\times \F_{q^2}^\times\mid x^3+y^3=1\}|\big/9\;\;\;,$$
since multiplying $x$ or $y$ by any cubic root of unity doesn't change $x^3$ nor $y^3$. It follows that $9l$ is almost equal to the number of points of the elliptic curve $x^3+y^3=z^3$ over $\F_{q^2}$: the difference consists of three points $(\theta,0,1)$ of the projective plane over $\F_{q^2}$, where $\theta$ is any cubic root of $1$, three points $(0,\theta,1)$, and three points $(\theta,-1,0)$. It follows that $9l=N_2-9$, where $N_2$ is the number of points over $\F_{q^2}$ of the Fermat cubic $E$ with equation $x^3+y^3=z^3$.\par
Now this is an elliptic curve, and by~\cite[(2.6)]{robert}, the zeta function of $E$ can be defined as
$$Z_E(u)={\rm exp}\big(\sum_{m\geq 1}N_m\fracb{u^m}{m}\big)\;\;\;,$$
where $N_m$ is the number of points of $E$ over $\F_{q^m}$. By~\cite[Theorem 2.8]{robert}, it has the following form
$$Z_E(u)=\fracb{1-au+qu^2}{(1-u)(1-qu)}\;\;\;,$$
 where $a=1+q-N_1$. Comparing the terms of degree 2 in $u$ in the expansion of those two expressions of $Z_E(u)$ as series in $u$ gives $N_2=N_1\big(2(q+1)-N_1\big)$.

Now since $3\mid q+1$, it follows that $3\nmid q-1$, and $x\mapsto x^3$ is a bijection of $\F_q$. Hence $E$ has as many points over $\F_q$ as the projective line with equation $x+y=z$, that is $N_1=q+1$. Hence $N_2=(q+1)^2$, which gives the following values for $l,n$ and $m$:
$$l=\Big(\fracb{q+1}{3}\Big)^2-1,\;\;\;m=\fracb{q^2-q-2}{9},\;\;\;n=
\Big(\fracb{q+1}{3}\Big)^2\;\;\;.$$
\end{Rem}

\begin{Theorem} \label{J2ZkN} Let $k$ be a field of characteristic $p$. Then:
\begin{enumerate}
\item The radical $J(ZkN)$ of the center of the group algebra $kN$ has a $k$-basis consisting of the elements $D_x$, for $x\in\F_{q^2}^\times/\Gamma-\{\Gamma\}$, $T_x$, for $x\in\Psi/\Gamma-\{\Gamma\}$, $T_1+\Id$, and $U_x$, for $x\in\F_{q^2}^\times/{\rm L}$. In particular, the dimension of $J(ZkN)$ is equal to $\fracb{q^2+q}{\gamma}+\gamma-1$.
\item The square $J^2(ZkN)$ of $J(ZkN)$ has a $k$ basis consisting of the elements $D_x+T_x$, where $x\in\Psi/\Gamma-\{\Gamma\}$. In particular, the dimension of $J^2(ZkN)$ is equal to $\fracb{q+1}{\gamma}-1$.
\item The cube $J^3(ZkN)$ of $J(ZkN)$ is equal to 0.
\end{enumerate}
\end{Theorem}

\begin{proof} As the group algebra $kN$ is indecomposable when $k$ is a field of characteristic $p$, the radical $J(ZkN)$ is equal to the kernel of the augmentation $\varepsilon:ZkN\to k$. If $X$ is the sum of the elements of a conjugacy class $C$ of $N$, then $\varepsilon(X)=|C|$, and by Proposition~\ref{tailles}, this is a multiple of $p$, unless $C$ is the class of the identity element of $N$, or $C$ is the class of $t_1$, and $|C|=q-1$ in this case. It follows that the elements listed in Assertion~(1) generate $J(ZkN)$. Moreover, they are obviously linearly independent, so they form a basis $\mathcal{B}$ of $J(ZkN)$.

Now by Proposition~\ref{DxT1}, for $x\in\Psi-\Gamma$, we have that $D_x(\Id+T_1)=D_x+T_x$ in $ZkN$, so the elements $D_x+T_x$, where $x\in\Psi/\Gamma-\{\Gamma\}$, are indeed in $J^2(ZkN)$, and they are clearly linearly independent. Moreover, reducing mod $p$ the formulas for products stated in Propositions~\ref{premiers produits} to~\ref{UxUy}, one checks easily that any product of two elements of the basis $\mathcal{B}$ is  equal to a (possibly zero) scalar multiple of an element $D_x+T_x$, for some $x\in\Psi-\Gamma$, and that the product of any three elements of $\mathcal{B}$ vanishes. This completes the proof of Theorem~\ref{J2ZkN}.
\end{proof}

If $k$ is a field of characteristic $p$ it is not difficult to give the explicit structure of $Z(kN)$
as a quotient of a polynomial ring in several variables.

\begin{Prop} Let $\gamma$ be the greatest common divisor of $3$ and $q+1$,
and let
$$\Gamma:=\{x\in\F_{q^2}\mid x^{\gamma}=1\}, \;\;\Psi:=\{x\in\F_{q^2}\mid x^{q+1}=1\},\;\;{\rm L}:=\{a^\gamma\mid a\in\F_{q^2}^\times\}\;\;\;.$$
Let ${\mathfrak U}:=\F_{q^2}^\times/\Gamma$, let ${\mathfrak V}:=\Psi/\Gamma$ and let ${\mathfrak W}:=\F_{q^2}^\times/{\rm L}$. Let $k$ be a field of characteristic $p>0$ and let $N$ be the normaliser of a Sylow $p$-subgroup of $PSU(3,q)$, where $p$ divides $q$.
Then, $$Z(kN)\simeq k[T,X_n,Y_m\;|\;n\in{\mathfrak W},m\in {\mathfrak U}]/I$$ where $I$ is the ideal generated by
$$T^2,\;TX_{n_1},\;TY_{m_1},\;X_{n_1}X_{n_2},\;X_{n_1}Y_{m_1},\;Y_{m_1}Y_{m_2},$$
$$X_{n_1}Y_{m_2}+\frac{1}{\gamma}X_{n_1}T,\;Y_{m_2}Y_{m_3}-(1-\delta_{m_2,m_3^{-1}})X_{m_2m_3}T$$
where $$n_1,n_2\in{\mathfrak W},m_1\in{\mathfrak U}-{\mathfrak V},m_2,m_3\in {\mathfrak V}$$
and $\delta_{a,b}$ is the Kronecker symbol.
\end{Prop}

\begin{proof}
We have a basis of $Z(kN)$ given in Theorem~\ref{J2ZkN} by the elements $D_x$, for $x\Gamma\in\F_{q^2}^\times/\Gamma$, $T_x$, for $x\Gamma\in\Psi/\Gamma-\{\Gamma\}$, $T_1+\Id$, and $U_x$, for $x{\rm L}\in\F_{q^2}^\times/{\rm L}$. Observe that $D_1=1$. Moreover, by Proposition~\ref{DxT1} we do not need to include $T_x$ as variable of the polynomial ring. This element is already the product of $T_1$ and $U_x$.

We obtain the following multiplication table.
$$\begin{array}{c||c|c|c|c|c|c}
&T_1+id&U_y&D_y\;(y\not\in\Psi)&D_y\;(y\in\Psi-\Gamma)\\
\hline\hline
T_1+id&0&0&0&T_y+D_y\\
&\text{\scriptsize Prop.~\ref{TxT1}}&\text{\scriptsize Prop.~\ref{T1Ux}}&\text{\scriptsize Prop.~\ref{premiers produits}(2)}&\text{\scriptsize Props.~\ref{DxT1}}\\\hline
U_x&0&0&0&-\frac{1}{\gamma}(D_x+T_x)\\
&\text{\scriptsize Prop.~\ref{T1Ux}}&\text{\scriptsize Prop.~\ref{UxUy}}&\text{\scriptsize Prop.~\ref{premiers produits}(3)}&\text{\scriptsize Prop.~\ref{DxUy}}\\\hline
D_x\;(x\not\in\Psi)&0&0&0&0\\
&\text{\scriptsize Prop.~\ref{premiers produits}(2)}&\text{\scriptsize Prop.~\ref{premiers produits}(3)}&\text{\scriptsize Props.~\ref{DxDyprodinPsi},\ref{premiers produits}(1),\ref{DxDxinv}}&\text{\scriptsize Prop.~\ref{premiers produits}(1)}\\\hline
D_x\;(x\in\Psi-\Gamma)&T_x+D_x&-\frac{1}{\gamma}(D_x+T_x)&0&(1-\delta_{xy\Gamma,\Gamma})(T_{xy}+D_{xy})\\
&\text{\scriptsize Prop.~\ref{DxT1}}&\text{\scriptsize Prop.~\ref{DxUy}}&\text{\scriptsize Prop.~\ref{premiers produits}(1)}&\text{\scriptsize Props.~\ref{DxDyPsi},\ref{DxDxinv}}\\
\end{array}$$
Now, mapping $T$ to $T_1+id$, $X_n$ to $U_n$ and $Y_m$ to $D_m$ gives an algebra homomorphism of the corresponding polynomial ring with kernel precisely the ideal $I$.
\end{proof}

\end{document}